\documentclass[12pt]{article}
\usepackage{amssymb, epsfig} 
\usepackage{amsmath} 
\usepackage{latexsym}

\newtheorem{Theorem}{Theorem}[section] 
\newtheorem{Definition}[Theorem]{Definition}

\newcommand{\R}{\ensuremath{\mathbb{R}}}

\newcommand{\1}{1\hspace{-1.4 mm}1} 
 
\begin{document}

                    \begin{center}
    {\large \bf WEAK SOLUTIONS FOR DISLOCATION TYPE    \\[0.1cm]
                                EQUATIONS }\\[1cm]
                    \end{center}

                    \begin{center}
                   {\sc Olivier Ley}\\
Laboratoire de Math\'ematiques et Physique Th\'eo\-ri\-que\\
CNRS UMR 6083, F\'ed\'eration Denis Poisson \\
Universit\'e de Tours,
Parc de Grandmont, 37200 Tours, France.\\
(ley@lmpt.univ-tours.fr)
\vspace{0.2cm}
                     \end{center}
\vspace{2cm}

\noindent
{\bf Abstract.} We describe recent results obtained by G. Barles, P. Car\-da\-lia\-guet, R. Monneau
and the author in \cite{bl06, bclm07}. They are concerned with
nonlocal Eikonal equations arising in the study of the 
dynamics of dislocation lines in crystals.
These equations are nonlocal but also non monotone.
We use a notion of weak solution to provide solutions
for all time. Then, we discuss the link between these weak solutions
and the classical viscosity solutions, and state some uniqueness
results in particular cases. 
A counter-example to uniqueness is given.

\vfill
\noindent
------------------------------------------------------------------------------
\\
{\footnotesize Communicated by xxxxxxxxxxxx; Received xxxxxxxxxx, 2007.\\
This work is supported by xxxxxxxxxxxxxxxxxxxx.
AMS Subject Classification 49L25, 35F25, 35A05, 35D05, 35B50, 45G10}


\section{Introduction} 
 
It is a great honor to contribute to this proceedings of the Conference
for the 25th Anniversary of Viscosity Solution and the Celebration of
the 60th birthday of Professor Hitoshi Ishii.

In this proceedings, we describe recent results \cite{bl06, bclm07}
obtained by the author in
collaboration with G. Barles, P. Cardaliaguet and R. Monneau
for first-order nonlocal Hamilton-Jacobi
modelling the dynamics of dislocations. 

Dislocations are defects in crystals
of typical length $10^{-6}m$ and the dynamics of dislocations is the main
microscopic explanation of the macroscopic behaviour of metallic crystals.
For details about the physics of dislocations, see for instance
Nabarro \cite{nabarro69} or Hirth and Lothe \cite{hl92}.
We are interested in a particular model introduced in Rodney, Le Bouar and Finel
\cite{rlf03}; the dislocation line evolves in a plane called
slip plane, with a normal velocity proportional to the 
Peach-Koehler force acting on this line. This Peach-Koehler force have two
contributions. The first one is the self-force created by the elastic field generated
by the dislocation line itself (i.e. this self-force is a nonlocal function of 
the shape of the dislocation line). The second one is due to exterior forces
(like an exterior stress applied on the material for instance).

More precisely, we study the evolution of a dislocation line $\Gamma_t$
which is, at any time $t\geq 0,$ the boundary of an open bounded set 
$\Omega_t\subset \R^N$ (with $N=2$ for the physical application).
The normal velocity, at each point $x\in\Gamma_t=\partial \Omega_t$
of the dislocation line, is given by
\begin{equation}\label{eq:1} 
V_n= c_0\star \1_{\overline\Omega_t} + c_1    
\end{equation} 
where $\1_{\overline\Omega_t}(x)$ is the indicator function of the set  
$\overline\Omega_t.$
The function $c_0(x,t)$ is a kernel  
which only depends on the physical properties 
of the crystal. 
In the special case of the study of dislocations, the kernel $c_0$ 
does not depend on time, but to keep a general setting we allow here a 
dependence on the time variable. 
Here $\star$ denotes the convolution in space, namely 
\begin{equation}\label{def-convolution}
(c_0(\cdot,t)\star \1_{\overline\Omega_t})(x)= \int_{\R^N} c_0(x-y,t)  
\1_{\overline\Omega_t}(y)dy, 
\end{equation} 
and this term appears to be the Peach-Koehler self-force created by the
dislocation itself, while $c_1(x,t)$ is the exterior contribution to the
velocity, created by
everything exterior to the dislocation line. We refer to Alvarez, Hoch, Le
Bouar and Monneau \cite{ahlm06} for a detailed presentation and a derivation of this 
model.

Using the level-set approach to front propagation problems, 
we can derive a partial differential equation 
to represent the evolution of $\Gamma_t.$ The level-set approach was
introduced by Osher and Sethian \cite{os88}, and then developped first by
Chen, Giga and Goto \cite{cgg91}, and Evans and Spruck \cite{es91}. This approach
produced a lot of applications and now there
is a huge literature; see the monograph of Giga \cite{giga06} for details.

The level-set approach consists in replacing the evolution of the  
set $\Gamma_t$ by the evolution of the zero level-set of an auxiliary 
function $u.$
More precisely, given a set $\Gamma_0$ (the dislocation line at time $t=0$)
and a bounded uniformly continuous function $u_0:\R^N\to \R$ such that 
\begin{eqnarray}\label{defu0}
\{u_0=0\} = \Gamma_0\ \ \ {\rm and} \ \ \ \{u_0>0\}=\Omega_0
\end{eqnarray}
($u_0$ represents the initial dislocation line), 
we are looking for a function $u : \R^N\times [0,T]\to \R$
which satisfies
\begin{eqnarray}\label{f-p}
\{u (\cdot ,t)=0\} = \Gamma_t\  \ {\rm and} \  \ \{u(\cdot ,t)>0\}=\Omega_t
\ \  \ {\rm for \ all \ } t\geq 0. 
\end{eqnarray}
The function $u$ has to satisfies the level-set equation 
(see \cite{giga06}) which reads here
\begin{equation}\label{dislocation} 
\left\{\begin{array}{l} 
\displaystyle{\frac{\partial u}{\partial t}   
= (c_0 (\cdot ,t) \star \1_{\{u(\cdot,t)\geq 0\}}(x) +c_1(x,t))|D u|  
\quad \mbox{in} 
\quad  \R^N\times (0,T)}\\ 
u(\cdot,0)=u_0\quad \mbox{in  }\R^N\; , 
\end{array}\right. 
\end{equation} 
where $\frac{\partial u}{\partial t},$ $Du$  and $|\cdot |$ denote respectively the
time and the spatial derivative of $u,$ and the Euclidean norm.  Note that
(\ref{def-convolution}) now reads 
\begin{eqnarray} \label{conv-space} 
c_0 (\cdot ,t) \star \1_{\{u(\cdot,t)\geq 0\}}(x)=\int_{\R^N} c_0(x-y,t)  
\1_{\{u(\cdot,t)\geq 0\}}(y)dy. 
\end{eqnarray} 
Note that (\ref{dislocation}) is not really a level-set equation because 
it is not invariant by increasing changes of functions. In order to have
rigorously a level-set equation, the
nonlocal term $\{u(\cdot,t)\geq 0\}$ should be replaced by
$\{u(\cdot,t)\geq u(x,t)\}$
(see Slep$\rm \check{c}$ev \cite{slepcev03}). But here, (\ref{dislocation})
is the equation we are interested in. 

The study of Equation (\ref{dislocation}) raises three main difficulties: the
first one is the presence of the nonlocal term (\ref{conv-space}).

The second difficulty is the weak regularity in time of the equation.
Indeed, as soon as $\{u (\cdot ,t)=0\}$ develops an interior 
(fattening phenomenon), the map
$t\mapsto c_0 (\cdot ,t) \star \1_{\{u(\cdot,t)\geq 0\}}(x)$ is no longer
continuous and we have to deal with (\ref{dislocation}) which
is an equation with measurable-in-time coefficients. 
The study of such equations was initiated by Ishii \cite{ishii85}
(see the Appendix).

The third difficulty, which is more involved, is a lack of monotonicity for~(\ref{dislocation}).
In many cases, proofs of existence and uniqueness for such
geometrical equations rely on the preservation of inclusion
property which can be stated as follows. 
Consider a front propagation problem (see (\ref{f-p})) with a given normal velocity.
Let $\Gamma_0$ and $\tilde{\Gamma}_0$ be two different initial fronts 
evolving independently. Then,
\begin{eqnarray}\label{pres-inclu}
\overline{\Omega}_0\subset {\rm int}(\tilde{\Omega}_0) \ \ \ 
\Longrightarrow \ \ \ 
\overline{\Omega}_t\subset {\rm int}(\tilde{\Omega}_t) \ 
{\rm for \ all \ time} \ t\geq 0.
\end{eqnarray}
Such a property is the key point to use the classical viscosity solutions'
theory. For instance, it is satisfied for local evolution problems as
propagation by constant normal velocities, mean curvature flow (see \cite{giga06})
or for some nonlocal problems as in Cardaliaguet \cite{cardaliaguet00, cardaliaguet01},
Dalio, Kim and Slep$\rm \check{c}$ev \cite{dlks04}, Srour \cite{srour06}, etc.
But, for dislocation dynamics, the kernel $c_0$ has a zero mean which implies that
it changes sign. Therefore, the preservation of inclusion property is not
true in general. It follows that we cannot expect a principle of comparison
(that is: the subsolutions of (\ref{dislocation}) are below the supersolutions).

For geometrical evolutions without preservation of inclusion, few results are
known, see however Giga, Goto and Ishii \cite{ggi92}, Soravia and Souganidis \cite{ss96} 
and Alibaud \cite{alibaud07}.
In the case of (\ref{dislocation}), under suitable assumptions  
on $c_0,c_1$ (see {\bf (H1)}-{\bf (H2)}) and on the initial data, the
existence and the uniqueness of 
the solution were proved first for short time in \cite{ahlm04, ahlm06}.
In \cite{acm05, bl06, cm06}, such results were proved
for all time under the additional assumption that $V_n\ge 0$, which is for 
instance always satisfied for $c_1$ satisfying 
$\displaystyle  c_1 (x,t) \geq |c_0(\cdot, t)|_{L^1(\R^N)}.$
 In the general case, a notion of weak solutions was introduced in
\cite{bclm07}.
 
The aim of this paper is to describe global-in-time results obtained in
 \cite{acm05, bl06, bclm07}.
In Section \ref{sec:existence12}, we define the weak solutions
and prove an existence theorem. In Section \ref{sec:uniq}, we state
some uniqueness results. Section \ref{sec:counter} is devoted
to the study of a counter-example to uniqueness.
Finally, we recall the definition of $L^1$-viscosity 
solutions and a new stability result proved by Barles \cite{barles06}
in the Appendix.

\section{Definition and existence of weak solutions}
\label{sec:existence12} 

We introduce the following notion of weak solutions for (\ref{dislocation}):

\begin{Definition}\label{defi:1}{\bf (Classical and weak solutions)} \cite{bclm07}\\ 
For any $T>0$, we say that a Lipschitz continuous function $u:\R^N\times
[0,T]\to \R$ 
is a weak solution of equation (\ref{dislocation}) on the time interval 
$[0,T)$, if there is some measurable map $\chi:\R^N\times (0,T)\to [0,1]$ such
that $u$ is a $L^1$-viscosity solution of 
\begin{equation}\label{FormeFaible} 
\left\{\begin{array}{l} 
\displaystyle{\frac{\partial u}{\partial t}  = {\bar c}(x,t) |D u| \quad \mbox{in} 
\quad  \R^N\times (0,T)}\\ 
u(\cdot,0)=u_0\quad \mbox{in  }\R^N\; , 
\end{array}\right. 
\end{equation} 
where  
\begin{equation}\label{eq:chi} 
{\bar c}(x,t) = c_0(\cdot,t)\star \chi(\cdot,t) (x) +c_1(x,t) 
\end{equation} 
and 
\begin{equation}\label{FormeFaible2} 
\begin{array}{l} 
\displaystyle{  \1_{\{u(\cdot,t) > 0\}}(x) \;  \leq \;  \chi(x,t) \; \leq  \1_{\{u(\cdot,t) \geq 0\}}(x) \; ,} 
\end{array} 
\end{equation} 
for almost all $(x,t)\in \R^N\times  [0,T]$. We say that $u$ is a classical solution of equation (\ref{dislocation}) if 
$u$ is a weak solution to (\ref{FormeFaible}) and if 
\begin{equation}\label{FormeFaible3} 
\1_{\{u(\cdot,t) > 0\}}(x) =\1_{\{u(\cdot,t) \geq 0\}}(x) 
\end{equation} 
for almost all $(x,t)\in \R^N\times  [0,T]$. 
\end{Definition} 

\noindent 
We recall that $L^1$-viscosity solutions were introduced by Ishii \cite{ishii85},
see the appendix for details.
Note that, for classical solutions, we have  $\chi(x,t)=\1_{\{u(\cdot,t) > 0\}}(x) =\1_{\{u(\cdot,t) \geq 0\}}(x)$ for 
almost all $(x,t)\in \R^N\times  [0,T].$ 
 
To state our first existence result, we introduce the following assumptions\\ 
 
\noindent {\bf (H0)} $u_0 : \R^N\to [-1,1]$ is Lipschitz continuous and
there exists $R_0>0$ such that $u_0(x) \equiv -1$  
for $|x|\geq R_0$, \\ 
 
\noindent {\bf (H1)} $c_0\in C([0,T); L^1\left(\R^N \right))$, $D_x c_0\in L^\infty ([0,T]; L^1\left(\R^N \right))$,  
$c_1\in C(\R^N\times [0,T])$ and there exists constants $M_1, L_1$ such that, for any $x,y\in \R^N$ and $t\in [0,T]$ 
\begin{eqnarray}\label{H1-1}
|c_1(x,t)|\le M_1 \quad \hbox{and} \quad |c_1(x,t)-c_1(y,t)|\le L_1|x-y| . 
 \end{eqnarray} 
Let us make some comments about these assumptions.
The role of $u_0$ is to represent the
initial dislocation $\Gamma_0$ which lies in a bounded region (see (\ref{defu0})).
In general, we choose $u_0$ as a truncation of the signed distance to $\Gamma_0$
(positive in $\Omega_0$). Such a function $u_0$ is Lispchitz continuous and
satisfies {\bf (H0)}.
Note that we do not impose any sign condition on $c_0$ in {\bf (H1)}.
In the sequel, we denote by $M_0, L_0$ some constants such that, for any (or almost every) $t\in [0,T)$, we have 
\begin{eqnarray}\label{H1-2}
|c_0(\cdot,t)|_{L^1(\R^N)}\le M_0 \quad \hbox{and} \quad 
|D_x c_0(\cdot ,t)|_{L^1(\R^N)} \le L_0.
 \end{eqnarray} 

Our first main result is the following. 
\begin{Theorem}\label{th:1}{\bf (Existence of weak solutions)} \cite{bclm07}\\ 
Under assumptions {\bf (H0)}-{\bf (H1)}, for any $T >0$ and for any initial data $u_0$,  
there exists a weak solution of equation (\ref{dislocation}) on the time 
interval $[0,T].$
\end{Theorem} 

We give only the main ideas of the proof of Theorem \ref{th:1}. The
whole proof can be found in \cite{bclm07} and an alternative proof is presented
in \cite{bclm07b}. \\

\noindent{\bf Sketch of proof of Theorem \ref{th:1}.} \\ 
\noindent{\it 1. Introduction of a perturbated equation.}  
We consider the equation 
\begin{eqnarray}\label{pert-eq} 
\hspace*{0.8cm} 
\frac{\partial u^\varepsilon}{\partial t}  = c_\varepsilon [u^\varepsilon](x,t) |D u^\varepsilon|  
\quad \mbox{in  }  \R^N\times (0,T)\; ,  
\end{eqnarray} 
where the unknown is $u^\varepsilon,$ 
\begin{eqnarray*}
c_\varepsilon [u]= \left( c_0 (\cdot,t)\star \psi_\varepsilon 
  (u (\cdot,t))  (x)+c_1(x,t)  
\right) \ \ \ {\rm for \ any} \ u:\R^N\times [0,T]\to \R,
\end{eqnarray*}
and $\psi_\varepsilon :\R\to \R$ is a sequence of continuous functions such  
that $\psi_\varepsilon (r) \equiv 0$ for $r \leq -\varepsilon$,  $\psi_\varepsilon  
(r) \equiv 1$ for $t \geq 0$  
and $\psi_\varepsilon$ is an affine function on $[-\varepsilon, 0]$. \\
 
\noindent{\it 2. Definition of a map $\mathcal{T}.$}  
Let
\begin{eqnarray*} 
X= \{ u\in C(\R^N\times [0,T]) : u\equiv -1 \ {\rm in} \ \R^N \backslash  
B(0,R_0 +MT), \\ 
\hspace*{4cm} |Du|, |\frac{\partial u}{\partial t}|/M \le |Du_0|_{L^\infty (\R^N)}  e^{L T}\},
\end{eqnarray*} 
where $M=M_0+M_1$ and $L=L_0+L_1$ (see (\ref{H1-1}) and (\ref{H1-2}) for the
definition of $M_0,M_1,L_0,L_1$). By Ascoli's Theorem, $X$ is a compact and
convex subset of $(C(\R^N\times [0,T]), |\cdot |_\infty).$
We define the map ${\mathcal T} : X\to X$ by~:  
if $u\in C(\R^N\times [0,T])$, then $u^\varepsilon:= {\mathcal T} (u)$ is the unique solution of  
(\ref{pert-eq}) with $c_\varepsilon [u]$ (instead of $c_\varepsilon
[u^\varepsilon]$). The existence and uniqueness of $u^\varepsilon$ come from 
classical results for Eikonal equations with finite speed propagation
property (see \cite[Theorem 2.1]{bclm07}, Crandall \& Lions \cite{cl83},  
\cite{ley01} and \cite{bl06}) since,  under assumption {\bf (H1)} 
on $c_1$ and $c_0$,  
$c_\varepsilon [u]$ satisfies {\bf (H1)} with fixed constants $M$ and $L.$   \\
 
\noindent{\it 3. Application of Schauder's fixed point theorem to  
$\mathcal{T}.$}  
The map ${\mathcal T}$ is continuous since $\psi_\varepsilon$ is continuous, by  
using the classical  
stability result for viscosity solutions (see Barles \cite{barles94}). 
Therefore, ${\mathcal T}$ has a fixed point $u_\varepsilon$ which is bounded in  
$W^{1,\infty} (\R^N\times [0,T])$  
uniformly with respect to $\varepsilon$ (since $M$ and $L$ are independent of  
$\varepsilon$). \\ 
 
\noindent{\it 4. Convergence of the fixed point when $\varepsilon\to 0.$}  
From Ascoli's Theorem, we extract a subsequence  
$(u_{\varepsilon'})_{\varepsilon'}$ which converges locally uniformly  
to a function denoted by $u.$ 
 The functions $\chi_{\varepsilon'}:=\psi_{\varepsilon'} (u_{\varepsilon'})$ 
satisfy $0\le\chi_{\varepsilon'} \le 1$. Therefore, we can extract a 
subsequence---still denoted $(\chi_{\varepsilon'})$---which  converges 
weakly$-*$ in $L^\infty_{\rm loc}(\R^N\times [0,T])$ to some function
$\chi:\R^N\times (0,T)\to [0,1]$.
Furthermore, setting $c_{\varepsilon'}= c_0 \star \chi_{\varepsilon'} + c_1,$ 
we have, for all $(x,t)\in \R^N\times [0,T],$ 
\begin{eqnarray*} 
\int_0^t c_{\varepsilon'}(x,s)ds  
&= & \int_0^t \int_{\R^N} c_0(x-y,s)\chi_{\varepsilon'} (y,s)dyds + \int_0^t c_1(x,s)ds \\ 
& \to & \int_0^t {\bar c}(x,s)ds, 
\end{eqnarray*} 
where ${\bar c}(x,t)=c_0(\cdot, t)\star \chi(\cdot,t) (x) + c_1(x,t).$ 
The above convergence is pointwise but, noticing that $c_{\varepsilon'}$ 
is bounded Lipschitz continuous in space uniformly in time and measurable
in time, we can apply the stability Theorem \ref{L1stability}
of Barles \cite{barles06} for
weak convergence in time. We conclude that  
$u$ is $L^1$-viscosity solution to (\ref{FormeFaible}) with ${\bar c}$ 
satisfying (\ref{eq:chi})-(\ref{FormeFaible2}). 
\hfill $\Box$ 

\section{Classical solutions and uniqueness results}
\label{sec:uniq}

Our second main result gives a sufficient condition for a weak solution to be a classical
one. 
\begin{Theorem}\label{th:2}{\bf (Links between weak solutions and 
classical continuous viscosity solutions)} \cite{bclm07}\\ 
Assume {\bf (H0)}-{\bf (H1)} and suppose that there is some $\delta\geq 0$ such that,
for all measurable map $\chi:\R^N\times (0,T)\to [0,1],$ 
\begin{equation}\label{HypPos} 
for \ all \ (x,t)\in \R^N \times [0,T],\quad  
c_0 (\cdot ,t) \star \chi (\cdot ,t) (x) +c_1(x,t)\geq \delta, 
\end{equation} 
and that the initial data $u_0$ satisfies (in the viscosity sense) 
\begin{eqnarray} \label{lgb1111} 
- |u_0| - |D u_0| \leq -\eta_0 \quad\hbox{in  } \R^N
\end{eqnarray} 
for some $\eta_0 >0$. Then any weak solution $u$ of 
(\ref{dislocation}) in the sense of Definition \ref{defi:1}, is a classical 
continuous viscosity solution of (\ref{dislocation}). 
\end{Theorem} 

Assumption (\ref{HypPos}) ensures that the velocity $V_n$ in (\ref{eq:1}) is 
nonnegative, i.e. the dislocation line is expanding. Of course, we can state similar 
results in the case of negative velocity for shrinking dislocation lines.
Assumption (\ref{lgb1111}) comes from  \cite{ley01}.
It means that $u_0$ is a viscosity 
subsolution of $-|v(x)|-|D v(x)|+\eta_0\leq 0.$  It can be seen as a 
nonsmooth generalization of the following situation:
if $u_0$ is $C^1,$ (\ref{lgb1111}) implies  that the gradient of $u_0$ does 
not vanish on the set $\{u_0=0\}$ and therefore this latter set is a $C^1$
hypersurface. \\

\noindent{\bf Sketch of proof of Theorem \ref{th:2}.}  
At first, if $u$ is a weak solution and ${\bar c}$ is 
associated with $u,$ then, from (\ref{HypPos}),
for any $x\in \R^N$ and for almost all $t\in [0,T]$, we have 
$$ 
 {\bar c}(x,t) \geq\delta \geq 0
$$
and therefore the Hamiltonian ${\bar c}(x,t)|Du|$ of (\ref{dislocation})
is convex $1$-homogeneous in the gradient variable.
Then, the conclusion 
is a consequence of a preservation of the lower-bound gradient
estimate (\ref{lgb1111}) proved in \cite[Theorem 4.2]{ley01}
for equations with convex Hamiltonians $H$ (such that
$H(x,t,\lambda p)=\lambda H(x,t,p)$ for all $\lambda \geq 0$):
there exists $\eta(T)>0$ such that
\begin{equation}\label{lbgT} 
- |u(\cdot ,t)| - |D u(\cdot ,t)| \leq -\eta(T) \quad\hbox{on  } \R^N\times (0,T). 
\end{equation} 
It follows that for every $t\in (0,T)$, the 
0--level-set of $u(\cdot,t)$ has a  zero Lebesgue measure
and therefore (\ref{FormeFaible3}) holds.
Moreover $t\mapsto \1_{\{u(\cdot,t) \geq 0\}}$ is also 
continuous in $L^1$, and then ${\bar{c}}$ is continuous.
\hfill $\Box$ \\

Let us turn to uniqueness results.
If the evolving set has positive velocity or if the velocity is nonnegative
and the following additional condition is fulfilled, then
 we can prove uniqueness results.\\ 
 
\noindent {\bf (H2)} $c_1$ and $c_0$ satisfy {\bf (H1)} and there exists 
constants $m_0, N_1$ and a positive function $N_0\in L^1(\R^N)$ such that, 
for any $ x,h\in \R^N$, $t\in [0,T)$, we have 
\begin{eqnarray*}
&& |c_0(x,t)| \le m_0,\\
&& c_1(x+h,t)+c_1(x-h,t)-2c_1(x,t)\ge -N_1 |h|^2,  \\
&& c_0(x+h,t)+ c_0(x-h,t)-2c_0(x,t) \ge -N_0(x) |h|^2.
\end{eqnarray*} 
Second and third conditions means that $c_0$ and $c_1$ are semiconvex in space. 

\smallskip 
 
\begin{Theorem}\label{th:3}{\bf (Uniqueness results)} \cite{acm05, bl06, bclm07}\\ 
Assume {\bf (H0)}-{\bf (H1)}-{\bf (H2)} and suppose that 
(\ref{HypPos}) and (\ref{lgb1111}) hold.
The solution of (\ref{dislocation})
is unique if 
 
(i) either $\delta=0$ and $u_0$ is semiconvex, i.e. satisfies for some constant $C>0$: 
$$
u_0(x+h)+u_0(x-h)-2u_0(x)  \ge -C |h|^2, \quad \forall x,h\in \R^N ;
$$ 

(ii) or $\delta>0.$
\end{Theorem} 

Even if it has no physical meaning in the theory of dislocations,
an important particular case of application of Theorem \ref{th:2} is the
uniqueness for (\ref{dislocation}) when $c_0\geq 0$ and $c_1\equiv 0$
(this implies (\ref{lgb1111})).
In this case, the preservation of inclusion property (\ref{pres-inclu})
holds and some classical results apply, see Cardaliaguet \cite{cardaliaguet00}
and \cite[Theorem 1.5]{bclm07}. But let us point out that
a nonnegative kernel $c_0$ does not ensure uniqueness in general,
see the counter-example in Section \ref{sec:counter}.

Point (i) of the theorem is the main result of \cite{acm05, bl06}. 
Let us compare the two articles.
In \cite{acm05}, it is proved that we have uniqueness for (\ref{dislocation})
if we start with an initial dislocation $\Gamma_0=\partial \Omega_0$ 
such that  $\Omega_0$ has the interior ball property of radius $r>0$ that is:
for any $x\in \overline{\Omega}_0,$ there exists $p\in\R^N\backslash\{0\}$
such that $\overline{B}(x-r\frac{p}{|p|},r)\subset \overline{\Omega}_0.$
In  \cite{bl06}, uniqueness is proved under the asumption that
$u_0$ is semiconvex and satisfies the lower-bound gradient (\ref{lgb1111}). 
This latter set of assumptions is equivalent to the interior ball property
for $\{u_0\geq 0\}$ (see \cite[Lemma A.1]{bl06}).\\

\noindent{\bf Sketch of proof of Theorem \ref{th:3}.} \\
\noindent{\it 1. Part (i) Definition of a map $\mathcal{F}.$} 
We follow the ideas of \cite{bl06} and refer to this paper for details. 
The proof relies on the Banach contraction fixed point theorem. Let
\begin{eqnarray*} 
Y&=& \{ \chi \in C([0,T], L^1(\R^N)) : \\
&& \hspace*{2.5cm} 0\leq \chi\leq 1,\, |\chi(\cdot ,t)|_{L^1(\R^N)}
\leq  \mathcal{L}^N (B(0,R_0 +MT)) \},
\end{eqnarray*} 
where $M=M_0+M_1$ (see (\ref{H1-1}) and (\ref{H1-2}) for the
definition of $M_0,M_1$), $\mathcal{L}^N$ is the Lebesgue measure in
$\R^N$ and $B(0,R)$ is the open ball of center $0$ and radius $R>0.$
For $\tau >0$ fixed, the set $Y$ is endowed with the norm
\begin{eqnarray*} 
|\chi|_{Y,\tau}= \mathop{\rm sup}_{t\in [0,\tau]}|\chi (\cdot ,t)|_{L^1(\R^N)}
\end{eqnarray*} 
Define $\mathcal{F}: Y\to Y$ by: for all $\chi\in Y,$
$\mathcal{F}(\chi)= \1_{u(\cdot ,t)\geq 0}$ where $u$ is the unique continuous
viscosity solution of 
\begin{equation}\label{dislo-gelee} 
\left\{\begin{array}{l} 
\displaystyle{\frac{\partial u}{\partial t}   
= c[\chi](x,t)|D u|  
\quad \mbox{in} 
\quad  \R^N\times (0,T)}\\ 
u(\cdot,0)=u_0\quad \mbox{in  }\R^N\; , 
\end{array}\right. 
\end{equation} 
where $c[\chi]= c_0 (\cdot ,t) \star \chi(\cdot, t)(x) +c_1(x,t).$
We have to check that $\mathcal{F}$ is well-defined. \\

\noindent{\it 2. Part (i) The map $\mathcal{F}$ is well defined.} 
From {\bf (H1)}-{\bf (H2)}, for all $\chi\in Y,$ 
the map $(x,t)\in \R^N\times [0,T]\mapsto c[\chi](x,t)$ is
bounded continuous in $\R^N\times [0,T],$ Lispchitz continuous and semiconvex in $x$ (uniformly
with respect to $t$) with some constants which depends only on the given
data $M_0,M_1,L_0,L_1,R_0, \eta(T), C.$ It follows that for all Lipschitz continuous $u_0,$
(\ref{dislo-gelee}) has a unique Lipschitz continuous viscosity solution
$u.$ Next,
if $u_0$ satisfies {\bf (H0)}, then, by the finite speed of propagation
property, for all $t\geq 0,$ $\{u(\cdot ,t)\geq 0\}\subset B(0,R_0+Mt).$
Let us give a geometrical interpretation of this latter property: by
(\ref{HypPos}), Equation (\ref{dislo-gelee}) is monotone and the preservation
inclusion principle (\ref{pres-inclu}) holds. Noticing that $M$ is an upper bound for 
the speed of propagation of the $0$-level-set of $u$ and  that $B(0,R_0+Mt)$
is the propagation of the ball $B(0,R_0)$ with normal velocity $M,$
by preservation of inclusion, the property follows.\\

\noindent{\it 3. Part (i) The map $\mathcal{F}$ is continuous.} 
It comes from 
the continuity of the map $t\in [0,T]\mapsto \int_{\R^N} \1_{u(\cdot ,t)\geq 0}(x)dx.$ 
The proof of this result is an immediate consequence of the preservation of the 
lower-bound gradient estimate (\ref{lgb1111})-(\ref{lbgT}) (see the proof of Theorem \ref{th:2}).\\

\noindent{\it 4. Part (i) Contraction property for $\mathcal{F}$ (beginning of
  the calculation).}
Let $\chi_1, \chi_2\in Y$ and $u_1, u_2$ be the solution of (\ref{dislo-gelee})
with $c[\chi_1]$ and $c[\chi_2]$ respectively. Set 
\begin{eqnarray}\label{defrho} 
\rho :=\mathop{\rm sup}_{t\in [0,\tau]} |(u_1-u_2)(\cdot ,t)|_{L^\infty (\R^N)}. 
\end{eqnarray} 
(note that $\rho\to 0$ as $\tau\to 0$ since $u_1(\cdot, 0)=u_2(\cdot,
0)=u_0$). For all $t\in [0,T],$ a straightforward computation leads to
\begin{eqnarray}\label{form986} 
\hspace*{1cm} &&|(\mathcal{F}(\chi_1)-\mathcal{F}(\chi_1))(\cdot ,t)|_{L^1 (\R^N)}\nonumber \\
&=&
|\1_{u_1(\cdot ,t)\geq 0}- \1_{u_2(\cdot ,t)\geq 0}|_{L^1 (\R^N)}\nonumber\\
&\leq & \mathcal{L}^N(\{ u_1(\cdot ,t)\geq 0, u_2(\cdot ,t)<0 \})
  + \mathcal{L}^N(\{ u_2(\cdot ,t)\geq 0, u_1(\cdot ,t)<0 \})\nonumber\\
&\leq &
\mathcal{L}^N(\{ -\rho \leq u_2(\cdot ,t)< 0\})
+\mathcal{L}^N(\{ -\rho \leq u_1(\cdot ,t)< 0\}).
\end{eqnarray}

\noindent{\it 5. Part (i) Contraction property for $\mathcal{F}$ ($L^1$-estimates).}
The estimate of the last two terms in (\ref{form986}) are based on some
fundamental $L^1$-estimates obtained in \cite{bl06}: 
let $\varphi_\varepsilon$ be a smooth approximation of $\1_{[-\rho,0)}$
(with $\1_{[-\rho,0)}\leq \varphi_\varepsilon\leq \1_{[-\rho-\varepsilon,\varepsilon]}$)
and $0< \rho < \eta(T)/2$ (where $\eta(T)$ is given by (\ref{lbgT})). Then,
there exists $K>0$ such that
\begin{eqnarray} \label{L1est}
\int_{\R^N} \varphi_\varepsilon(u_2(x,t))dx 
\leq {\rm e}^{Kt} \int_{\R^N} \varphi_\varepsilon(u_0(x)dx 
\end{eqnarray} 
which implies by sending $\varepsilon\to 0,$
\begin{eqnarray*} 
\mathcal{L}^N(\{ -\rho \leq u_2(\cdot ,t)< 0\})
\leq {\rm e}^{Kt}\mathcal{L}^N(\{ -\rho \leq u_0< 0\})
\end{eqnarray*} 
(we have the same formula for $u_1$).
We provide a formal calculation which emphasizes the main ideas
(see \cite[Proposition 3.1]{bl06} for a rigorous computation).
We have
\begin{eqnarray*}
\frac{d}{dt}\left( \int_{\R^N} \varphi_\varepsilon (u_2(x,t))dx\right)
= \int_{\R^N} \varphi_\varepsilon ' (u_2(x,t)) \frac{\partial u_2}{\partial t}(x,t) dx
\end{eqnarray*}
for a.e. $t\in [0,T].$
Using Equation (\ref{dislo-gelee}), it follows
\begin{eqnarray*}
\int_{\R^N} \varphi_\varepsilon ' (u_2) \frac{\partial u_2}{\partial t} dx
& = &
\int_{\R^N} \varphi_\varepsilon ' (u_2) c[\chi_2]|D u_2|dx \\
&=&
\int_{\R^N} \langle \varphi_\varepsilon ' (u_2) D u_2 , \frac{c[\chi_2]D u_2}{|D u_2|}\rangle dx \\
&=&
\int_{\R^N} \langle D \varphi_\varepsilon  (u_2), \frac{c[\chi_2]D u_2}{|D u_2|}\rangle dx
\end{eqnarray*}
since, from $0< \rho < \eta(T)/2,$ and (\ref{lbgT}), we have
$|D u_2|>\eta(T)/2,$ for almost every $(x,t)$ such that 
$\varphi  (u_2(x,t))\not=0.$ By an integration by parts, we obtain
\begin{eqnarray*}  
\hspace*{1cm}
\int_{\R^N} \langle D \varphi_\varepsilon  (u_2), \frac{c[\chi_2]D u_2}{|D u_2|}\rangle dx
=
- \int_{\R^N} \varphi_\varepsilon  (u_2) \, {\rm div}( c[\chi_2] \frac{D u_2}{|D u_2|}) dx.
\end{eqnarray*}
Since $u_0$ is semiconvex and {\bf (H1)}-{\bf (H2)} hold, by \cite[Theorem 5.2]{ley01},
the solutions $u_1,u_2$ of (\ref{dislo-gelee}) are still semiconvex in space, i.e.
there exists $\overline{C}$ such that
\begin{eqnarray*}
D^2 u_1, D^2 u_2 \geq -\overline{C}\, Id \ \ {\rm for \ a.e.} \
(x,t)\in\R^N\times [0,T]. 
\end{eqnarray*}
Using this estimate and the lower-bound gradient estimate (\ref{lbgT}) again, we have,
for almost every $(x,t)\in \R^N\times [0,T]$ such that $\varphi  (u_2(x,t))\not=0,$
\begin{equation*} 
  -  {\rm div}( \frac{D u_2}{|D u_2|})
= -\frac{1}{|D u_2|} {\rm trace} \left[ \left( Id-\frac{D u_2\otimes D u_2}
{|D u_2|^2}\right) D^2 u_2\right]
\leq \frac{2 \overline{C}}{\eta(T)}.
\end{equation*}
It gives
\begin{eqnarray*}  
-  {\rm div}( c[\chi_2] \frac{D u_2}{|D u_2|})
&=& -\langle D c[\chi_2], \frac{D u_2}{|D u_2|}\rangle
- c[\chi_2]\,  {\rm div}( \frac{D u_2}{|D u_2|}) \nonumber \\
&\leq &
|D c[\chi_2]|_{L^\infty(\R^N)} + \frac{2 \overline{C}|c[\chi_2]|_{L^\infty(\R^N)}}{\eta(T)}:=K
\end{eqnarray*}
since $c[\chi_2]$ is nonnegative bounded Lipschitz continuous by
(\ref{HypPos}) and Step 1.
Finally, we obtain, for a.e. $t\in[0,T],$
\begin{eqnarray*}
\frac{d}{dt}\left( \int_{\R^N} \varphi_\varepsilon (u_2(x,t))dx\right)
\leq K \int_{\R^N} \varphi_\varepsilon (u_2(x,t))dx
\end{eqnarray*}
which yields (\ref{L1est}) through a classical Gronwall's argument.

By the same kind of arguments, we can estimate
$\mathcal{L}^N(\{ -\rho \leq u_0< 0\})$ to obtain
\begin{equation}\label{somln}
\mathcal{L}^N(\{ -\rho \!\leq \!u_2(\cdot ,t)\!< \!0\})
+\mathcal{L}^N(\{ -\rho \!\leq \!u_1(\cdot ,t)\! <\!  0\})
\leq \frac{2C}{\eta_0} \mathcal{L}^N(B(0,R_0\! +\! 1)){\rm e}^{Kt}\, \rho.
\end{equation}

\noindent{\it 6. Part (i) Contraction property for $\mathcal{F}$ (stability
estimates with respect to variations of the velocity).}
Since $u_1$ and  $u_2$ are the solutions of (\ref{dislo-gelee})
with $c[\chi_1]$ and $c[\chi_2]$ respectively, we have the ``continuous
dependence'' type result: for all $t\in [0,T],$
\begin{equation}\label{contdep}
|(u_1-u_2)(\cdot ,t)|_{L^\infty (\R^N)}\leq |Du_0|_{L^\infty (\R^N)} {\rm
  e}^{\Lambda t} \int_0^t |(c[\chi_1]-c[\chi_2])(\cdot ,s)|_{L^\infty (\R^N)} ds,
\end{equation}
where $\Lambda= {\rm max}\{ |D c[\chi_1]|_{L^\infty (\R^N)}, |D c[\chi_2]|_{L^\infty (\R^N)}\}.$

\noindent{\it 7. Part (i) Contraction property for $\mathcal{F}$ (end of the proof).}
From (\ref{defrho}), (\ref{form986}),  (\ref{somln}) and (\ref{contdep}), we get
\begin{eqnarray*}
&& |\mathcal{F}(\chi_1)-\mathcal{F}(\chi_1)|_{Y,\tau}\\
&\leq& \frac{2C}{\eta_0} \mathcal{L}^N(B(0,R_0\! +\! 1)){\rm e}^{Kt} 
\mathop{\rm sup}_{t\in [0,\tau]} |(u_1-u_2)(\cdot ,t)|_{L^\infty (\R^N)}\\
&\leq &
\overline{L} \mathop{\rm sup}_{t\in [0,\tau]}\int_0^t |(c[\chi_1]-c[\chi_2])(\cdot ,s)|_{L^\infty (\R^N)} ds\\
&\leq&\overline{L}\tau |\chi_1- \chi_2|_{Y,\tau}
\end{eqnarray*}
for some constant $\overline{L}.$
Therefore, we have contraction for $\tau$ small enough. This implies the 
uniqueness of a classical solution to (\ref{dislocation}) on the time interval
$[0,\tau].$
Noticing that all the constants depend only on the given data, we conclude by
a step-by-step argument to obtain the uniqueness on the whole interval $[0,T].$ \\

\noindent{\it 8. Part (ii).} 
The additional difficulty comparing to the proof of (i) is the fact that
$u_0$ is not supposed to be semiconvex anymore and then $u(\cdot, t)$ is not
semiconvex.
Nevertheless, we assume that $\delta >0,$ i.e. the velocity is
positive. Such a property implies the creation of the interior ball property of radius $\gamma
t$ for $\{ u(\cdot, t)\geq 0\}$ for every $t>0$ (see Cannarsa and Frankowska \cite{cf06}
and \cite[Lemma 2.3]{bclm07}). Roughly speaking, we recover this way the semiconvexity
property for $u(\cdot ,t)$ (see the comment after the statement of Theorem
\ref{th:3}).

Using arguments similar to those in the proof of
Part (i) and the interior ball regularization,
we prove the following Gronwall type inequality
\begin{eqnarray*}
&& | \1_{\{u_1(\cdot,t) \geq 0\}}- \1_{\{u_2(\cdot,t) \geq 0\}}|_{L^1(\R^N)}\\ 
&\leq & \displaystyle{   C \left[  {\rm per}(\{u_1(\cdot,t)\geq 0\})
+{\rm per}(\{u_2(\cdot,t)\geq 0\})\right] } \\
&& \hspace*{4cm}
\int_0^t | \1_{\{u_1(\cdot,s) \geq 0\}} - \1_{\{u_2(\cdot,s) \geq 0\}}|_{L^1(\R^N)}ds  
\end{eqnarray*} 
where $u_i,$ $i=1,2$ are two weak solutions of (\ref{dislocation}),
$C$ is a constant depending on the constants of the problem and ${\rm per}(\{u_i(\cdot,t)\geq 0\})$ is 
the $\mathcal{H}^{N-1}$ measure (the perimeter) of  the set $\partial \{u_i(\cdot,t)\geq 0\}).$ 
In order to apply Gronwall's Lemma it is sufficient to know that 
the functions $t\mapsto  {\rm per}(\{u_i(\cdot,t)\geq 0\})$ 
belong to $L^1$. This fact is proved by applying the co-area formula.
Finally, it follows $\1_{\{u_1(\cdot,t)\geq 0\}}=\1_{\{u_2(\cdot,t)\geq 0\}}$ for all $t\in [0,T]$
and therefore $u_1=u_2$ since they are solution of the same equation.
\hfill $\Box$ 

\section{A counter-example to uniqueness \cite{bclm07}}
\label{sec:counter}

The following example is inspired from \cite{bss93}.

Let us consider, in dimension $N=1$, the following equation of type 
(\ref{dislocation}), 
\begin{equation}\label{dislocation-d1} 
\left\{\begin{array}{l} 
\displaystyle{\frac{\partial U}{\partial t}   
= (1 \star \1_{\{U(\cdot,t)\geq 0\}}(x) +c_1(t))|D U|  
\quad \mbox{in} 
\quad  \R\times (0,2]}\\[2mm] 
U(\cdot,0)=u_0\quad \mbox{in  }\R\; , 
\end{array}\right. 
\end{equation} 
where we set 
$c_0(x,t):= 1,$ $c_1(x,t):=c_1(t)= 2(t-1)(2-t)$ and $u_0(x)=1-|x|.$ 
Note that $1 \star \1_A ={\mathcal L}^1 (A)$ for any measurable 
set $A\subset \R.$ 

Note that $c_0\equiv 1$ 
does not satisfies exactly {\bf (H1)} but this is not the point
here: because of the finite speed of
propagation property, it is possible to modify $c_0$ such that
{\bf (H1)} and the construction below holds.
 
We start by solving auxiliary problems for time in $[0,1]$ and $[1,2]$ in order to 
produce a family of solutions for the original problem in $[0,2].$ 
 
\noindent 
{\it 1. Construction of a solution for $0\leq t\leq 1.$} The function $x_1(t)=(t-1)^2$ is the solution of 
the ordinary differential equation (ode in short) 
\begin{eqnarray*} 
\dot{x}_1 (t)= c_1(t)+2x_1(t) \ {\rm for} \ 0\leq t\leq 1, \ \ \ 
{\rm and } \ x(0)=1, 
\end{eqnarray*}  
(note that $\dot{x}_1\leq 0$ in $[0,1]$). 
Consider 
\begin{eqnarray} \label{eq-eik-1} 
\left\{\begin{array}{l} 
\displaystyle{\frac{\partial u}{\partial t}= \dot{x}_1 (t)\left|\frac{\partial 
  u}{\partial x}\right| \quad\hbox{in  } \R \times (0,1],} \\ 
u(\cdot,0)=u_0\quad\hbox{in  } \R. 
\end{array}\right. 
\end{eqnarray}  
There exists a unique continuous viscosity solution $u$ of (\ref{eq-eik-1}). 
Looking for $u$ under the form $u(x,t)=v(x,\Gamma (t))$ with $\Gamma (0)=0,$ 
we obtain that $v$ satisfies 
\begin{eqnarray*} 
\frac{\partial v}{\partial t}\, \dot{\Gamma}(t)= \dot{x}_1 (t)\left|\frac{\partial 
  v}{\partial x}\right|. 
\end{eqnarray*}  
Choosing $\Gamma (t)=-x_1(t)+1,$ we get that $v$ is the solution of 
\begin{eqnarray*}  
\left\{\begin{array}{l} 
\displaystyle{\frac{\partial v}{\partial t}= -\left|\frac{\partial 
  v}{\partial x}\right| \quad\hbox{in  } \R \times (0,1],} \\ 
v(\cdot,0)=u_0\quad\hbox{in  } \R. 
\end{array}\right. 
\end{eqnarray*}  
By the Oleinik-Lax formula, $\displaystyle v(x,t)=\mathop{\rm inf}_{|x-y|\leq t} 
u_0(y).$ Since $u_0$ is even, we have, for all $(x,t)\in\R\times [0,1],$ 
\begin{eqnarray*} 
u(x,t)= \mathop{\rm inf}_{|x-y|\leq \Gamma (t)} u_0(y) =u_0 (|x|+\Gamma (t))= 
u_0 (|x|-x_1(t)+1). 
\end{eqnarray*}  
Therefore, for $0\leq t\leq 1,$ 
\begin{equation} \label{incl678} 
\{ u(\cdot ,t)>0\}= (-x_1(t),x_1(t)) \  \ \ {\rm and} \  \ \ 
\{ u(\cdot ,t)\geq 0\}= [-x_1(t),x_1(t)]. 
\end{equation}  
We will see in Step 3 that $u$ is a solution of (\ref{dislocation-d1}) in $[0,1].$ 
 
\noindent 
{\it 2. Construction of solutions for $1\leq t\leq 2.$} Consider now, 
for any measurable function $0\leq \gamma (t)  \leq 1,$ the unique 
solution $y_\gamma$ of the ode 
\begin{eqnarray}\label{edo123} 
\dot{y}_\gamma (t)= c_1(t)+2\gamma (t)  y_\gamma (t) \ {\rm for} \ 1\leq t\leq 2, \ \ \ 
{\rm and } \  y_\gamma (1)=0. 
\end{eqnarray}  
By comparison, we have $0\leq y_0(t) \leq y_\gamma (t)\leq y_1 (t)$ for $1\leq 
t\leq 2,$ where $y_0, y_1$ are the solutions of (\ref{edo123}) obtained with $\gamma (t)\equiv 0,1.$ 
In particular, it follows that $\dot{y}_\gamma \geq 0$ in $[1,2].$ 
Consider 
\begin{eqnarray*}  
\left\{\begin{array}{l} 
\displaystyle{\frac{\partial u_\gamma}{\partial t}= \dot{y}_\gamma (t)\left|\frac{\partial 
  u_\gamma}{\partial x}\right| \quad\hbox{in  } \R \times (1,2],} \\ 
u_\gamma(\cdot,1)=u (\cdot ,1)\quad\hbox{in  } \R, 
\end{array}\right. 
\end{eqnarray*}  
where $u$ is the solution of (\ref{eq-eik-1}). 
Again, this problem has a unique continuous viscosity solution $u_\gamma$ and 
setting 
$\Gamma_\gamma (t)=y_\gamma (t)\geq 0$ for $1\leq t\leq 2,$ we obtain that 
$v_\gamma$ defined by $v_\gamma (x,\Gamma_\gamma (t))=u_\gamma (x,t)$ is the 
unique continuous viscosity solution of 
\begin{eqnarray*}  
\left\{\begin{array}{l} 
\displaystyle{\frac{\partial v_\gamma}{\partial t}= \left|\frac{\partial 
  v_\gamma}{\partial x}\right| \quad\hbox{in  } \R \times (0,\Gamma_\gamma (2)],} \\ 
v_\gamma(\cdot,0)=u (\cdot ,1)\quad\hbox{in  } \R. 
\end{array}\right. 
\end{eqnarray*}  
Therefore, for all $(x,t)\in\R\times [1,2],$ we have 
\begin{eqnarray*} 
u_\gamma(x,t)= \mathop{\rm sup}_{|x-y|\leq y_\gamma (t)} u(y,1) = 
\left\{\begin{array}{ll} 
0 & {\rm if} \ |x|\leq y_\gamma (t), \\ 
u(|x|-y_\gamma (t), 1) & {\rm otherwise}. 
\end{array}\right. 
\end{eqnarray*}  
(Note that $u(-x,t)=u(x,t)$ since $u_0$ is even and, since $u(\cdot ,1)\leq 
0,$ by the maximum principle, we have $u_\gamma \leq 0$ in $\R\times [1,2].$) 
It follows that, for all $1\leq t\leq 2,$ 
\begin{equation} \label{incl679} 
\{ u_\gamma(\cdot ,t)>0\}= \emptyset \  \ \ {\rm and} \  \ \ 
\{ u_\gamma(\cdot ,t)\geq 0\}= \{ u_\gamma (\cdot ,t)=0\} = [-y_\gamma (t),y_\gamma (t)]. 
\end{equation}  
 
\begin{figure}[ht] 
\begin{center} 
\epsfig{file=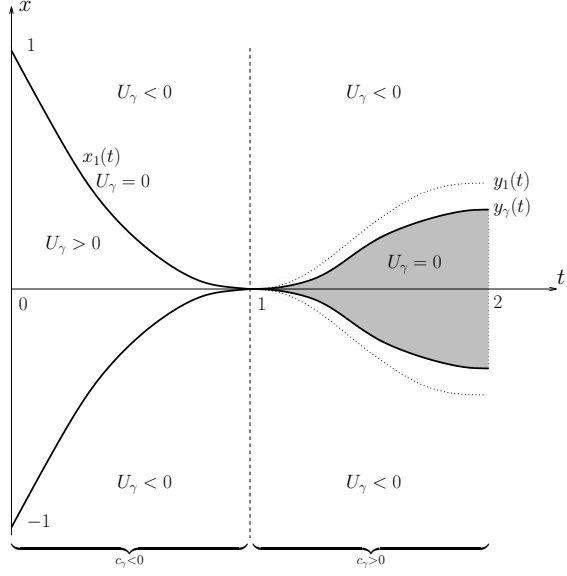, width=7cm}  
\end{center} 
\caption{\label{dess-gamma} 
{\it Fattening phenomenon for the functions $U_\gamma.$}} 
\end{figure}

\noindent 
{\it 3. There are several weak solutions of (\ref{dislocation-d1}).} 
Set, for $0\leq \gamma (t)  \leq 1,$ 
\begin{eqnarray*} 
\begin{array}{lll} 
c_\gamma (t)= c_1(t)+2 x_1(t), & U_\gamma (x,t)= u(x,t) &  {\rm if} \ (x,t)\in \R\times [0,1], \\  
c_\gamma (t)= c_1(t)+2 \gamma (t) y_\gamma (t), & U_\gamma (x,t)= u_\gamma (x,t) &  {\rm if} \ (x,t)\in \R\times [1,2]. \\ 
\end{array} 
\end{eqnarray*}  
Then, from Steps 1 and 2,  $U_\gamma$ is the unique continuous viscosity 
solution of  
\begin{eqnarray*} 
\left\{\begin{array}{l} 
\displaystyle{\frac{\partial U_\gamma}{\partial t}=  c_\gamma (t) \left|\frac{\partial 
  U_\gamma}{\partial x}\right| \quad\hbox{in  } \R \times (0,2],} \\ 
U_\gamma(\cdot,0)=u_0 \quad\hbox{in  } \R. 
\end{array}\right. 
\end{eqnarray*}  
Taking $\chi_\gamma (\cdot ,t)=\gamma (t) \1_{[ -y_{\gamma}(t), y_{\gamma} (t) ]}$ for $1\leq t\leq 2,$ 
from (\ref{incl678}) and (\ref{incl679}), we have 
\begin{eqnarray*} 
\1_{\{U_\gamma (\cdot ,t)>0\}}\leq \chi_\gamma (\cdot ,t) \leq \1_{\{U_\gamma (\cdot ,t)\geq 0\}}, 
\end{eqnarray*}  
(see Figure \ref{dess-gamma}). 
It follows that all the $U_\gamma$'s, for measurable $0\leq \gamma (t) \leq 1,$ are weak 
solutions of (\ref{dislocation-d1}) so we do not have uniqueness and the set of solutions is quite large.  

\section*{Appendix:  $L^1$-viscosity solutions and 
a stability result for weak convergence in time}   
\label{appendice} 

We recall that the definition of $L^1$-viscosity solutions was introduced 
in Ishii's paper \cite{ishii85}. We refer also to Nunziante \cite{nunziante90, nunziante92} 
and Bourgoing \cite{bourgoing04a,bourgoing04b} for a complete presentation of the theory. 
 
Consider the equation 
\begin{equation}\label{eq:cbar} 
\left\{\begin{array}{l} 
\displaystyle{\frac{\partial v}{\partial t}  = {\bar{c}}(x,t)|D v| \quad \mbox{in} 
\quad  \R^N\times (0,T)}\\ 
\\ 
v(\cdot,0)=u_0\quad \mbox{in  }\R^N\; , 
\end{array}\right. 
\end{equation} 
where the velocity ${\bar{c}} : \R^N\times (0,T) \to \R$ is defined for 
almost every $t\in (0,T)$. We 
also assume that ${\bar{c}}$ satisfies\\ 
\noindent {\bf (H3)} The function ${\bar{c}}$ is continuous with respect to $x\in \R^N$ 
and measurable in $t.$  
For all $x,y\in\R^N$ and almost all $t\in [0,T],$  
$$ 
|{\bar{c}}(x,t)|\le M \ \ \ {\rm and} \ \ \ 
|{\bar{c}}(x,t)-{\bar{c}}(y,t)|\le L|x-y|. 
$$ 

\begin{Definition}{\bf ($L^1$-viscosity solutions)}\\ 
An upper-semicontinuous (respectively lower-semicontinuous) function $v$  on $\R^N\times [0,T]$  
is a $L^1$-viscosity subsolution (respectively supersolution) of 
(\ref{eq:cbar}), if  
$$ 
v(0,\cdot) \le u_0 \quad (\mbox{respectively}\quad v(0,\cdot) \ge u_0), 
$$ 
and if for every $(x_0,t_0)\in \R^N\times [0,T]$, $b\in L^1(0,T)$, $\varphi \in 
C^\infty(\R^N\times (0,T))$ and continuous function $G: \R^N\times (0,T) \times \R^N\to \R$ 
such that \\ 
(i) the function  
$$ 
(x,t)\longmapsto v(x,t)-\int_0^t b(s)ds -\varphi(x,t) 
$$  
has a local maximum 
(respectively minimum) at $(x_0,t_0)$ over $\R^N\times (0,T)$ and such that \\ 
(ii) for almost every $t\in (0,T)$ in some neighborhood of $t_0$ and for 
every $(x,p)$ in some neighborhood of $(x_0,p_0)$ with $p_0=D\varphi(x_0,t_0)$, 
we have 
$${\bar{c}}(x,t)|p| -b(t)  \le G(x,t,p) \quad (\mbox{respectively}\quad 
{\bar{c}}(x,t)|p| -b(t)  \ge G(x,t,p))$$ 
then 
$$ 
\frac{\partial \varphi}{\partial t}(x_0,t_0) \le G(x_0,t_0,p_0)\quad 
(\mbox{respectively}\quad  \frac{\partial \varphi}{\partial t}(x_0,t_0) \ge 
G(x_0,t_0,p_0)). 
$$ 
Finally we say that a locally bounded function $v$ defined on $\R^N\times [0,T]$ 
is a $L^1$-viscosity solution of (\ref{eq:cbar}), if its upper-semicontinuous 
(respectively lower-semicontinuous) envelope is a 
$L^1$-viscosity subsolution (respectively supersolution). 
\end{Definition} 

 
\begin{Theorem}{\bf (Existence and uniqueness in the $L^1$ sense)}\\ 
For any $T>0$, under assumptions {\bf (H0)} and {\bf (H3)}, there exists a 
unique $L^1$-viscosity solution to (\ref{eq:cbar}). 
\end{Theorem} 
 
Finally, let us consider the solutions $v^\varepsilon$ to the following 
equation 
\begin{equation}\label{eq:cbarbis} 
\left\{\begin{array}{l} 
\displaystyle{\frac{\partial v^\varepsilon}{\partial t}  = {\bar{c}}^\varepsilon(x,t)|D v^\varepsilon| \quad \mbox{in} 
\  \R^N\times (0,T)},\\[2mm] 
v^\varepsilon(\cdot,0)=u_0\quad \mbox{in}\ \R^N\;. 
\end{array}\right. 
\end{equation} 
The following stability result is a particular case of a 
general stability result proved by Barles in \cite{barles06}. 
 
\begin{Theorem}\label{L1stability}{\bf ($L^1$-stability)} \cite{barles06}\\ 
Under assumption {\bf (H0)}, let us assume that the velocity 
${\bar{c}}^\varepsilon$ satisfies {\bf (H3)} (with some constants $M,L$ 
independent of $\varepsilon$). Let us consider the $L^1$-viscosity solution $v^\varepsilon$ to (\ref{eq:cbarbis}). Assume that 
$v^\varepsilon$ converges locally uniformly to a function $v$ and, for all $x\in \R^N,$ 
\begin{eqnarray*}
\hspace*{0.5cm} \int_0^t {\bar{c}}^\varepsilon (x,s)ds \to \int_0^t {\bar{c}} (x,s)ds  
\ \ \ { locally \ uniformly \ in } \ (0,T). 
\end{eqnarray*} 
Then $v$ is a $L^1$-viscosity solution of (\ref{eq:cbar}). 
\end{Theorem} 



\begin{thebibliography}{10}

\bibitem{alibaud07}
N.~Alibaud.
\newblock Existence, uniqueness and regularity for nonlinear degenerate
  parabolic equations with nonlocal terms.
\newblock {\em To appear in NoDEA Nonlinear Differential Equations Appl.}

\bibitem{acm05}
O.~Alvarez, P.~Cardaliaguet, and R.~Monneau.
\newblock Existence and uniqueness for dislocation dynamics with nonnegative
  velocity.
\newblock {\em Interfaces Free Bound.}, 7:415--434, 2005.

\bibitem{ahlm04}
O.~Alvarez, P.~Hoch, Y.~Le~Bouar, and R.~Monneau.
\newblock R\'esolution en temps court d'une \'equation de {H}amilton-{J}acobi
  non locale d\'ecrivant la dynamique d'une dislocation.
\newblock {\em C. R. Math. Acad. Sci. Paris}, 338(9):679--684, 2004.

\bibitem{ahlm06}
O.~Alvarez, P.~Hoch, Y.~Le~Bouar, and R.~Monneau.
\newblock Dislocation dynamics: short-time existence and uniqueness of the
  solution.
\newblock {\em Arch. Ration. Mech. Anal.}, 181(3):449--504, 2006.

\bibitem{barles94}
G.~Barles.
\newblock {\em Solutions de viscosit\'e des \'equations de
  {H}amilton-{J}acobi}.
\newblock Springer-Verlag, Paris, 1994.

\bibitem{barles06}
G.~Barles.
\newblock A new stability result for viscosity solutions of nonlinear parabolic
  equations with weak convergence in time.
\newblock {\em C. R. Math. Acad. Sci. Paris}, 343(3):173--178, 2006.

\bibitem{bclm07}
G.~Barles, P.~Cardaliaguet, O.~Ley, and R.~Monneau.
\newblock Global existence results and uniqueness for dislocation equations.
\newblock {\em To appear in SIAM J. Math. Anal.}

\bibitem{bclm07b}
G.~Barles, P.~Cardaliaguet, O.~Ley, and A.~Monteillet.
\newblock {\em In preparation}, 2007.

\bibitem{bl06}
G.~Barles and O.~Ley.
\newblock Nonlocal first-order {H}amilton-{J}acobi equations modelling
  dislocations dynamics.
\newblock {\em Comm. Partial Differential Equations}, 31(8):1191--1208, 2006.

\bibitem{bss93}
G.~Barles, H.~M. Soner, and P.~E. Souganidis.
\newblock Front propagation and phase field theory.
\newblock {\em SIAM J. Control Optim.}, 31(2):439--469, 1993.

\bibitem{bourgoing04a}
M.~Bourgoing.
\newblock Vicosity solutions of fully nonlinear second order parabolic
  equations with ${L}^1$-time dependence and {N}eumann boundary conditions.
\newblock {\em To appear in NoDEA Nonlinear Differential Equations Appl.}

\bibitem{bourgoing04b}
M.~Bourgoing.
\newblock Vicosity solutions of fully nonlinear second order parabolic
  equations with ${L}^1$-time dependence and {N}eumann boundary conditions.
  existence and applications to the level-set approach.
\newblock {\em To appear in NoDEA Nonlinear Differential Equations Appl.}

\bibitem{cf06}
P.~Cannarsa and H.~Frankowska.
\newblock Interior sphere property of attainable sets and time optimal control
  problems.
\newblock {\em ESAIM Control Optim. Calc. Var.}, 12(2):350--370 (electronic),
  2006.

\bibitem{cardaliaguet00}
P.~Cardaliaguet.
\newblock On front propagation problems with nonlocal terms.
\newblock {\em Adv. Differential Equations}, 5(1-3):213--268, 2000.

\bibitem{cardaliaguet01}
P.~Cardaliaguet.
\newblock Front propagation problems with nonlocal terms. {II}.
\newblock {\em J. Math. Anal. Appl.}, 260(2):572--601, 2001.

\bibitem{cm06}
P.~Cardaliaguet and C.~Marchi.
\newblock Regularity of the eikonal equation with {N}eumann boundary conditions
  in the plane: application to fronts with nonlocal terms.
\newblock {\em SIAM J. Control Optim.}, 45(3):1017--1038 (electronic), 2006.

\bibitem{cgg91}
Y.~G. Chen, Y.~Giga, and S.~Goto.
\newblock Uniqueness and existence of viscosity solutions of generalized mean
  curvature flow equations.
\newblock {\em J. Differential Geom.}, 33(3):749--786, 1991.

\bibitem{cl83}
M.~G. Crandall and P.-L. Lions.
\newblock Viscosity solutions of {H}amilton-{J}acobi equations.
\newblock {\em Trans. Amer. Math. Soc.}, 277(1):1--42, 1983.

\bibitem{dlks04}
F.~Da~Lio, C.~I. Kim, and D.~Slep{\v{c}}ev.
\newblock Nonlocal front propagation problems in bounded domains with
  {N}eumann-type boundary conditions and applications.
\newblock {\em Asymptot. Anal.}, 37(3-4):257--292, 2004.

\bibitem{es91}
L.~C. Evans and J.~Spruck.
\newblock Motion of level sets by mean curvature. {I}.
\newblock {\em J. Differential Geom.}, 33(3):635--681, 1991.

\bibitem{giga06}
Y.~Giga.
\newblock {\em Surface evolution equations}, volume~99 of {\em Monographs in
  Mathematics}.
\newblock Birkh\"auser Verlag, Basel, 2006.
\newblock A level set approach.

\bibitem{ggi92}
Y.~Giga, S.~Goto, and H.~Ishii.
\newblock Global existence of weak solutions for interface equations coupled
  with diffusion equations.
\newblock {\em SIAM J. Math. Anal.}, 23(4):821--835, 1992.

\bibitem{hl92}
J.~R. Hirth and L.~Lothe.
\newblock {\em Theory of dislocations}.
\newblock Krieger, Malabar, Florida, second edition, 1992.

\bibitem{ishii85}
H.~Ishii.
\newblock {H}amilton-{J}acobi equations with discontinuous {H}amiltonians on
  arbitrary open sets.
\newblock {\em Bull. Fac. Sci. Eng. Chuo Univ.}, 28:33--77, 1985.

\bibitem{ley01}
O.~Ley.
\newblock Lower-bound gradient estimates for first-order {H}amilton-{J}acobi
  equations and applications to the regularity of propagating fronts.
\newblock {\em Adv. Differential Equations}, 6(5):547--576, 2001.

\bibitem{nabarro69}
F.~R.~N. Nabarro.
\newblock {\em Theory of crystal dislocations}.
\newblock Clarendon Press, Oxford, 1969.

\bibitem{nunziante90}
D.~Nunziante.
\newblock Uniqueness of viscosity solutions of fully nonlinear second order
  parabolic equations with discontinuous time-dependence.
\newblock {\em Differential Integral Equations}, 3(1):77--91, 1990.

\bibitem{nunziante92}
D.~Nunziante.
\newblock Existence and uniqueness of unbounded viscosity solutions of
  parabolic equations with discontinuous time-dependence.
\newblock {\em Nonlinear Anal.}, 18(11):1033--1062, 1992.

\bibitem{os88}
S.~Osher and J.~Sethian.
\newblock Fronts propagating with curvature dependent speed: algorithms based
  on {H}amilton-{J}acobi formulations.
\newblock {\em J. Comp. Physics}, 79:12--49, 1988.

\bibitem{rlf03}
D.~Rodney, Y.~Le~Bouar, and A.~Finel.
\newblock Phase field methods and dislocations.
\newblock {\em Acta Materialia}, 51:17--30, 2003.

\bibitem{slepcev03}
D.~Slep{\v{c}}ev.
\newblock Approximation schemes for propagation of fronts with nonlocal
  velocities and {N}eumann boundary conditions.
\newblock {\em Nonlinear Anal.}, 52(1):79--115, 2003.

\bibitem{ss96}
P.~Soravia and P.~E. Souganidis.
\newblock Phase-field theory for {F}itz{H}ugh-{N}agumo-type systems.
\newblock {\em SIAM J. Math. Anal.}, 27(5):1341--1359, 1996.

\bibitem{srour06}
A.~Srour.
\newblock Nonlocal second-order {H}amilton-{J}acobi equations arising in
  tomographic reconstruction.
\newblock {\em submitted}, 2006.

\end{thebibliography}
\end{document}